\documentstyle{article}
\begin{document}
\baselineskip+8pt
\small \begin{center}{\textit{ In the name of
 Allah, the Beneficent, the Merciful.}}\end{center}
\begin{center}{\bf  A particular case of the Jacobian conjecture}
\end{center}
\begin{center}{\bf  Ural Bekbaev  \\Turin Polytechnic University in Tashkent,\\
INSPEM, Universiti Putra Malaysia.\\
e-mail: bekbaev2011@gmail.com, u.bekbaev@polito.uz}
 \end{center}
\begin{abstract}{A particular case of the Jacobian conjecture is considered and for small dimensional cases a computational approach is offered.   }\end{abstract}


\vspace{0.5cm}

In this paper we want to show an application of a result from [1] and offer a computational approach, in small dimensional cases, to the Jacobian conjecture.

Let $F$ stand for the field of real or complex numbers,  $\varphi: F^n \rightarrow F^n$ and $f: F^n \rightarrow F^{n'}$ be any polynomial maps. In such case one
 can attach to $f$ the following sequence of polynomial maps $\Phi_m : F^n \rightarrow F^{n'}$ defined  by the recurrent formula
$\Phi_{m+1}(f)=\Phi_m(f)-\Phi_m(f)\circ\varphi$, where $\Phi_0(f) =f$ and $\circ$ stands for the composition(superposition)
operation.
It is not difficult to see that \ \ \ $$\Phi_m(f)=\sum_{k=0}^{m}(-1)^{k}\left(\begin{array}{c}
  m \\
  k\\
\end{array}\right)f\circ\varphi^{\circ k}$$

{\bf Proposition 1.} \emph{The following properties of $\Phi_m$ are evident:}

1. $\Phi_m(\lambda f)=\lambda \Phi_m(f)$ \emph{for any} $\lambda\in F$

2. $\Phi_m(f+g)= \Phi_m(f)+\Phi_m(g)$, \emph{where} $g: F^n \rightarrow F^{n'}$

3. $\Phi_m(f)\circ \varphi=\Phi_m(f\circ \varphi)$

4. $\Phi_m(\Phi_k(f))=\Phi_{m+k}(f)$

In [1] the following criterion of polynomiality was given.

{\bf Proposition 2.} \emph{If $\varphi: F^n \rightarrow F^n$ is a polynomial map of the form $\varphi(x)=x+\mbox{"higher order terms"}$ then $\varphi^{-1}(x)$ is a polynomial map whenever $\Phi_m(id)=0$ for some natural $m$, where $id: F^n \rightarrow F^n$ is the identity map $id(x)=x$, $x=(x_1,x_2,...,x_n)$.}

In common case checking polynomiality of $\varphi^{-1}(x)$ by finding explicit expression for $\varphi^{-1}(x)$ in terms of $x$ is not an easy thing. In such cases the the presented criterion may be very helpful.

We don't know if polynomiality of $\varphi(x)$ and $\varphi^{-1}(x)$ implies existence of such a natural $m$ for which $\Phi_m(id)=0$ but if one wants to find a counter example to the Jacobian Conjecture [2], probably such example exists, he can try to find it among such polynomials $\varphi(x)=x+\mbox{"higher order terms"}$ for which the Jacobian determinant is 1 but for it never $\Phi_m(id)=0$. Moreover according to [2]  to prove the Jacobian Conjecture it is
enough to prove it for polynomial maps $\varphi(x): F^n
\rightarrow F^n$ of the form $\varphi(x)=x+((xA_1)^3,
(xA_2)^3,..., (xA_n)^3)$, where $x=(x_1,x_2,...,x_n)$-row vector of independent variables,
$A_i$ are the corresponding columns of a matrix
$A=(a^i_j)_{i,j=1,2,...,n}$. In this paper we follow this direction.

{\bf Definition 1.} \emph{A polynomial map $f: F^n \rightarrow F$ is said to be
$\varphi$-invariant if $f\circ \varphi=f$ that is $\Phi_1(f)=0$.}

Let $I_0(\varphi)$ stand for the set of all such $\varphi$-invariant polynomials. It is clear that $\Phi_m(\lambda f)=\lambda \Phi_m(f)$ for any $\lambda(x)\in I_0(\varphi)$

{\bf Definition 2.} \emph{A polynomial map $f: F^n \rightarrow F$ is said to be
nearly $\varphi$-invariant if $f\circ \varphi\in f+I_0(\varphi)$ that is $\Phi_2(f)=0$.}

Let $I_1(\varphi)$ stand for the set of all such nearly $\varphi$-invariant polynomials.

{\bf Proposition 3.} \emph{If $f\in I_1(\varphi)$ then for any natural $m\geq 1$ one has
$\Phi_{m+1}(f^m)=0$.}

{\bf Proof.} If $f\circ \varphi= f+f_0$, where $f_0\in I_0(\varphi)$, then $f^m-(f\circ \varphi)^m=\sum_{i=0}^{m-1}f^ic_i$, where $c_i\in I_0(\varphi)$.
Let us prove Proposition 3 by induction on $m$. Proposition 3 is true at $m=1$ due to the condition $\Phi_{1+1}(f^1)=0$ and $\Phi_{m+1}(f^m)=\Phi_{m}(f^m)-\Phi_{m}(f^m)\circ \varphi=
\Phi_{m}(f^m-f^m\circ \varphi)=\Phi_{m}(f^m-(f\circ \varphi)^m)=\Phi_{m}(\sum_{i=0}^{m-1}f^ic_i)=\sum_{i=0}^{m-1}c_i\Phi_{m}(f^i)=0$ on induction.

The following result is about $n=2$  and $n=3$ cases.

{\bf Theorem.} \emph{If $\varphi(x)=x+((xA_1)^3, (xA_2)^3)$ or $\varphi(x)=x+((xA_1)^3, (xA_2)^3,(xA_3)^3)$ and determinant of its Jacobian matrix equals 1 then $\varphi^{-1}(x)$ is a polynomial in  $x$.}

{\bf Proof.} Here is a proof of $n=3$ case, $n=2$ case is more simple. The equality of
determinant of its Jacobian matrix to 1 yields in the following
system of equalities.

\begin{equation} \left\{\begin{array}{l}
 (xA_1)^2a_1^1+ (xA_2)^2a_2^2+(xA_3)^2a_3^3=0 \\
 (xA_1)^2(xA_2)^2\left|
\begin{array}{cc}
  a_1^1 & a_2^1 \\
  a_1^2 & a_2^2 \\
\end{array}
\right|+ (xA_1)^2(xA_3)^2\left|
\begin{array}{cc}
  a_1^1 & a_3^1 \\
  a_1^3 & a_3^3 \\
\end{array}
\right|+(xA_2)^2(xA_3)^2\left|
\begin{array}{cc}
  a_2^2 & a_3^2 \\
  a_2^3 & a_3^3 \\
\end{array}
\right|=0 \\
(xA_1)^2(xA_2)^2(xA_3)^2\left|
\begin{array}{ccc}
  a_1^1 & a_2^1& a_3^1\\
  a_1^2 & a_2^2&a_3^2 \\
  a_1^3 & a_2^3&a_3^3 \\
\end{array}\right|=0\\
\end{array}\right.\end{equation}

Let $r$ stand for the rank of the matrix $A$

1-case. If $r=1$ then, for example $A_1\neq 0$ and
$A_2=\lambda_2A_1,A_3=\lambda_3A_1$, system of equalities (1)
is equivalent to
$$a_1^1+ \lambda_2^2a_1^2+\lambda_3^2a_1^3=0$$ , $\varphi(x)=x+(xA_1)^3(1, \lambda_2^3, \lambda_3^3)$, $\Phi_{1}(id)(x)=-(xA_1)^3(1, \lambda_2^3, \lambda_3^3)$
and $f(x)=xA_1$ is a $\varphi$ -invariant. It implies that $\Phi_{2}(id)=0$ and therefore $\varphi^{-1}$ is a polynomial map.

2-case. If $r=2$ then, for example $\{A_1, A_2\}$ is linear
independent and $A_3=\lambda_1A_1+\lambda_2A_2$, system (1) is
equivalent to
\begin{equation} \left\{\begin{array}{l}
 a_1^1+ \lambda_1^2a_3^3= a_2^2+ \lambda_2^2a_3^3= \lambda_1\lambda_2a_3^3=0 \\
 \left|
\begin{array}{cc}
  a_1^1 & a_2^1 \\
  a_1^2 & a_2^2 \\
\end{array}
\right|+\lambda_1^2\left|\begin{array}{cc}
  a_2^2 & a_3^2 \\
  a_2^3 & a_3^3 \\
\end{array}
\right|+ \lambda_2^2\left|
\begin{array}{cc}
  a_1^1 & a_3^1 \\
  a_1^3 & a_3^3 \\
\end{array}
\right|=0\\
\lambda_2^2\left|\begin{array}{cc}
  a_2^2 & a_3^2 \\
  a_2^3 & a_3^3 \\
\end{array}
\right|=\lambda_1^2\left|
\begin{array}{cc}
  a_1^1 & a_3^1 \\
  a_1^3 & a_3^3 \\
\end{array}
\right|=\lambda_1\lambda_2\left|
\begin{array}{cc}
  a_1^1 & a_3^1 \\
  a_1^3 & a_3^3 \\
\end{array}
\right|=\lambda_1\lambda_2\left|
\begin{array}{cc}
  a_2^2 & a_3^2 \\
  a_2^3 & a_3^3 \\
\end{array}
\right|=0 \\
 \end{array}\right.\end{equation}

2-1-case. If $\lambda_1=0$ then (2) is equivalent to
$a_1^1=a_2^2+\lambda_2^3a_2^3=a_2^1(a_1^2+\lambda_2^3a_1^3)=0$.

2-1-1-case. If in addition to 2-1-case $a_2^1=0$ then
$$\varphi_1(x)=x_1+(a_1^2x_2+a_1^3x_3)^3,\varphi_2(x)=x_2+(a_2^3(-\lambda_2^3x_2+x_3))^3,
\varphi_3(x)=x_3+(\lambda_2a_2^3(-\lambda_2^3x_2+x_3))^3$$
 This map is of the form
$\varphi_1(x)=x_1+(ax_2+bx_3)^3,\varphi_2(x)=x_2+c(-\lambda x_2+x_3)^3,
\varphi_3(x)=x_3+c\lambda(-\lambda x_2+x_3)^3$ for which $\Phi_{1}(id)(x)=-((ax_2+bx_3)^3, c(-\lambda x_2+x_3)^3, c\lambda(-\lambda x_2+x_3)^3)$,
 $-\lambda x_2+x_3\in I_0(\varphi)$ and
$f(x)=ax_2+bx_3\in I_1(\varphi)$. It implies that at least $\Phi_{5}(id)=0$ for this $\varphi$ and therefore $\varphi^{-1}$ is a polynomial map.

2-1-2-case. If in addition to 2-1-case $a_2^1\neq 0$ then
$$\varphi_1(x)=x_1+(a_1^3(-\lambda_2^3x_2+x_3))^3,\varphi_2(x)=x_2+(x_1a_2^1+a_2^3(-\lambda_2^3x_2+x_3))^3,
\varphi_3(x)=x_3+\lambda_2^3(x_1a_2^1+a_2^3(-\lambda_2^3x_2+x_3))^3$$
 This map is of the form
$$\varphi_1(x)=x_1+a(-\lambda x_2+x_3)^3,\varphi_2(x)=x_2+(bx_1+c(-\lambda x_2+x_3))^3,
\varphi_3(x)=x_3+\lambda(bx_1+c(-\lambda x_2+x_3))^3$$
  Once again $-\lambda x_2+x_3\in I_0(\varphi)$ and
   $f(x)=bx_1\in I_1(\varphi)$, so at least $\Phi_{5}(id)=0$ for this $\varphi$ and therefore $\varphi^{-1}$ is a polynomial map.

2-2-case. If $\lambda_1\lambda_2\neq 0$ then (2) is equivalent to $a_3^3=a_1^1=a_2^2=a_1^2a_2^1=a_1^3a_3^1=a_3^2a_2^3=0$ and there is no such kind polynomial map $\varphi$.

\begin{center}{References}\end{center}

[1] Ural Bekbaev. \emph{An inversion formula for multivariate power series}.\\ arXiv:1203.3834v1 [math. AG] 17 Mar 2012

[2] Dan Yan, Michiel de Bondt. \emph{Some remarks on the Jacobian conjecture and polynomial endomorphisms}. \\ arXiv:1203.3609v1 [math. AG] 16 Mar 2012

\end{document}